\documentstyle[12pt,leqno,amssymb,amscd,graphics,hyperref]{amsart}  
\hypersetup{hidelinks=true}

\newtheorem{thm}{Theorem}[section]

\newtheorem{prop}[thm]{Proposition}

\newtheorem{clar}[thm]{Clarification}
\theoremstyle{definition}

\theoremstyle{remark}

\begin{document}

\newcommand{\ct}{\cite}
\newcommand{\pr}{\protect\ref}
\newcommand{\su}{\subseteq}
\newcommand{\pa}{{\partial}}
\newcommand{\e}{\epsilon}
\newcommand{\es}{{\varnothing}}

\newcommand{\D}{{\lambda}}

\newcommand{\Q}{{\mathbb Q}}
\newcommand{\R}{{\mathbb R}}
\newcommand{\rn}{{{\mathbb R}^n}}
\newcommand{\Z}{{\mathbb Z}}
\newcommand{\N}{{\mathbb N}}

\newcommand{\h}{\hat}
\newcommand{\dd}{\delta}
\newcommand{\x}{\Delta}
\newcommand{\xk}{{\x^k_{v_1,\dots,v_k}}}
\newcommand{\xt}{{\x^2_{v,w}}}
\newcommand{\sm}{\s  \! M}
\newcommand{\sx}{\s  \! X}

\newcommand{\1}{^{1 \over \D}_{\phantom i}}

\newcommand{\s}{{}^*}

\newcommand{\f}{{\varphi}}

\newcommand{\ha}{{\mathfrak{h}}}
\newcommand{\cp}{{\mathcal{P}}}
\newcommand{\cu}{{\mathcal{U}}}
\newcommand{\cd}{{\mathcal{D}}}
\newcommand{\cs}{{\mathcal{S}}}
\newcommand{\cc}{{\mathcal{C}}}
\newcommand{\lf}{{\lfloor}}
\newcommand{\rf}{{\rfloor}}
\newcommand{\vv}{\overrightarrow}
\newcommand{\tf}{{\widetilde{F}}}
\newcommand{\tg}{{\widetilde{G}}}
\newcommand{\der}{{\partial}}

\newcommand{\th}{{   {}^{\mathfrak{h}}     }}

\newcommand{\tm}{{{}^{\mathfrak{h}} \! {M}}}

\newcommand{\ty}{{{}^{\mathfrak{h}}  {Y}}}
\newcommand{\tu}{{{}^{\mathfrak{h}}  {U}}}
\newcommand{\tr}{{{}^{\mathfrak{h}}  {\R}}}
\newcommand{\ii}{{(0,1)}}

\makeatletter\def\l@subsection{\@tocline{2}{0pt}{4pc}{5pc}{}}\makeatother

\newcounter{numb}

\title{Resolution of the Surprise Exam Paradox}

\author{Tahl Nowik}
\address{Department of Mathematics, Bar-Ilan University, 
Ramat-Gan 5290002, Israel}
\email{tahl@@math.biu.ac.il}
\urladdr{\url{www.math.biu.ac.il/~tahl}}

\date{November 30, 2014}
\begin{abstract}
We present a resolution of the celebrated ``Surprise Exam Paradox''. We argue that 
if the surprise exam story  is analyzed using the exact same meaning of the notion of ``surprise'' as is dictated by the story itself, then no paradox arises.
\end{abstract}

\maketitle

The ``Surprise Exam Paradox'' has drawn great interest  over the years. An overview  appears in \ct{c}, with extensive bibliography which is being updated in its arXiv version. I first describe the surprise exam story, as a play in three acts.

\section{The surprise exam story}\label{intro}

{\bf{ACT 1:}} Teacher in class: ``On one of the days of next week, I will give you a surprise exam!''

{\bf{ACT 2:}} Student A: ``If the exam is not given before Friday, then Friday morning before school, we would know that the exam must be on that day, and so it would not be a surprise. Since the teacher has announced that the exam will be a surprise, he will not give it on Friday. 

Student B: ``I now imagine Thursday morning.  
Since we have already established that the exam will not be on Friday, if the exam is not given before Thursday, then Thursday morning before school, we would know that the exam must be on that day, and so it would not be a surprise. Since the teacher has announced that the exam will be a surprise, he will not give it on Thursday either. 

Narrator: ``... and so the students go on, analyzing the previous day, and the one preceding that, until the first day of the week, and are thus convinced that the teacher's surprise exam cannot take place on any day of the week...''

{\bf{ACT 3:}}  On Wednesday the teacher announces: ``Today we have an exam!''

 Narrator: ``... and the students were so very surprised...''

\section{Apparent paradox, and its resolution}

The seeming paradox that this story suggests, is that though the students have established that a surprise exam cannot take place, they are surprised on Wednesday after all, as the narrator announces in Act 3. 
I however will argue, that the narrator's notion of ``surprise'' appearing in Act 3 is different from the students' notion of surprise appearing in Act 2 when deducing that no surprise exam can take place. 
I will establish that if the narrator would have held the 
precise same meaning of the notion of surprise  as used by the students in Act 2,  then in Act 3 he would say: 
``and as expected, the students were not surprised!",
and so, no paradox would arise.

We ask then, how is the notion of surprise used by the students in Act 2 in their deduction that no surprise exam can take place. In the students' argument in Act 2, the notion of surprise is used in exactly one crucial step of the inductive argument. When for a given day they realize that in the morning of that day, before school, they will be able to prove to themselves that the exam must take place on that day, then \emph{they 
feel that an exam under such circumstance is not a surprise}, and so they deduce that the exam cannot take place on that day, and another step of the induction is completed.
So we see that in fact the notion ``\emph{not} being surprised'' or ``unsurprised'' 
 appears, and it is used precisely as follows:

\begin{prop}\label{p}
The notion of ``not being surprised'' is used by the students in Act 2 in the following way:
\emph{If in the morning of the exam, before school, we can provide ourselves with a proof that the exam must take place on that day, then an exam on that day is not a surprise.} 
\end{prop}

At first sight this seems as an acceptable meaning for the notion of unsurprised,
and we do not question it when following the students' arguments in Act 2. I will however argue, that this meaning of unsurprised is quite different from
our natural notion which is  represented  by the narrator in Act 3.
To establish this claim, I will formulate some clarifications of Proposition \pr{p}.

\begin{clar}\label{c1} In the eyes of the students the following holds: \emph{
Even if in addition to the proof that the exam \emph{must} take place today, we can also provide ourselves with a proof that the exam \emph{cannot} take place today, we still consider an exam given  today \emph{not} to be a surprise.}
\end{clar}

How do we know this? We know this because we see it happening in Act 2, though it is slightly hidden by 
differences of time in the story. The students imagine themselves in the future, on Wednesday morning before school, formulating the proof that the exam \emph{must} take place on that day, and then, using their notion of unsurprised as appears in Proposition \pr{p}, and the teacher's announcement of Act 1, 
they infer in the present time of Act 2 that the exam \emph{cannot} take place on Wednesday. Since we see that our students are able to imagine themselves making the first deduction (``exam must take place'') on Wednesday, they can also imagine themselves making their second, contradicting  deduction (``exam cannot take place'') on Wednesday. 
In fact, when they imagine Tuesday morning, which is one step ahead in the induction, they indeed imagine themselves realizing that there cannot be an exam on Wednesday.
So, there is no reason why they shouldn't imagine themselves making that same deduction on Wednesday.
We see however that this does not deter them from feeling that an exam on Wednesday would not be a surprise, since we see  they continue untroubled with their inductive argument. This establishes the claim of Clarification \pr{c1}.

We  go further in demonstrating the extent to which our students hold  the notion of unsurprised as formulated in
Proposition \pr{p}. 

\begin{clar}\label{c2}  In the students eyes the following holds: \emph{
Even if  in several previous mornings we have provided ourselves with a proof that the exam must take place on that given day, and then have seen that it has in fact not taken place, we are not discouraged. That is,  if today we may again provide ourselves with a proof that the exam must take place today, 
we still feel that an exam taking place today is not a surprise.}
\end{clar}

And how do we know this? In the same way that we have established Clarification \pr{c1}.
During Act 2, the students can imagine themselves  proving to themselves every morning until Tuesday that the exam must take place on that given day, which is something they in fact do during Act 2 (in reverse order), and then imagine the exam not taking place on each of these days, which is also a fact that they happen to establish in Act 2. And so they know that when Wednesday morning arrives, and they will formulate their proof that the exam must take place on that day, they will have already experienced several failures of their predictions.
We see again that this does not deter them from feeling that an exam on that day is not  a surprise, and the induction process proceeds. This establishes the claim of Clarification \pr{c2}.

We can go on and point to other peculiarities of the notion of unsurprised appearing in Proposition 
\pr{p}, and then formulate a corresponding clarification that says that even if this peculiarity presents itself, the students continue to hold their precise view of unsurprised appearing in Proposition \pr{p}.
The argument will be similar to the argument appearing twice above. Whatever peculiarity we note, can also be noted by the students. Since we see they have continued confidently with their inductive process, we see
that in their eyes, any such peculiarity does not invalidate their notion of unsurprised appearing in Proposition \pr{p}. 

Now, in view of our understanding of the students'  notion of unsurprised, we can move on to analyze Act 3. If the narrator would use precisely
the same notion of unsurprised as held by the students, he would announce in Act 3 that the students were not
surprised, since their notion of unsurprised does indeed  hold when the teacher announces the exam on Wednesday. On Wednesday morning they can indeed provide themselves with a proof that the exam must take place that day (the proof described in Act 2), and so according to their particular notion of unsurprised formulated in Proposition \pr{p}, and further clarified in Clarifications \pr{c1} and \pr{c2}, they will indeed \emph{not} be surprised.

In conclusion, the source of the apparent paradox  is a crucial difference between our intuitive notion of surprise, and that of the students appearing in Act 2. 
We are puzzled by this story since when following the arguments of Act 2 we are willing to accept the 
notion of surprise used by the students, not noticing its peculiarities, and so we also accept their conclusion, that no surprise exam can take place.
But then in Act 3, we slip back, together with the narrator, to some different, more common, presumably less peculiar notion of surprise,
according to which an exam on Wednesday \emph{is} a surprise. 
Since these are two different notions of surprise, their contradicting conclusions are not a paradox.


\begin{thebibliography}{StoPC}



\bibitem[C]{c}
T. W. Chow: ``The Surprise Examination or Unexpected Hanging Paradox''
\emph{Amer. Math. Monthly} 105 (1998), 41--51. Also posted at arXiv:math/9903160

\end{thebibliography}
\end{document}